\documentclass{gtart_h}

\def\ifplaintex{\expandafter\ifx\csname documentclass\endcsname\relax}

\def\gtp{{\mathsurround=0pt\it $\cal G\mskip-2mu$eometry \&\ 
$\cal T\!\!$opology $\cal P\!$ublications}}  

\def\recd{{\small Received:\qua\receiveddate\ifx\reviseddate\relax
\else\qquad Revised:\qua\reviseddate\fi\par}} 

\def\lognumber#1{\def\thelognumber{#1}}
\def\volumenumber#1{\def\thevolumenumber{#1}}
\def\volumeyear#1{\def\thevolumeyear{#1}}
\def\papernumber#1{\def\thepapernumber{#1}}
\def\pagenumbers#1#2{\def\startpage{#1}\def\finishpage{#2}}
\def\published#1{\def\publishdate{#1}}

\def\received#1{\def\receiveddate{#1}}

\def\accepted#1{\def\accepteddate{#1}}
\def\asciititle#1{\def\theasciititle{#1}}
\def\covertitle#1{\def\thecovertitle{#1}}

\def\asciiaddress#1{\def\theasciiaddress{#1}}
\def\asciiemail#1{\def\theasciiemail{#1}}

\long\def\asciiabstract#1{\long\def\theasciiabstract{#1}}
\def\asciikeywords#1{\def\theasciikeywords{#1}}

\let\\\par\let\thelognumber\relax\let\thevolumenumber\relax
\let\thepapernumber\relax\let\thevolumeyear\relax\let\startpage\relax
\let\finishpage\relax\let\publishdate\relax\let\receiveddate\relax
\let\reviseddate\relax\let\accepteddate\relax\let\theasciititle\relax
\let\thecovertitle\relax\let\theasciiauthors\relax\let\theasciiaddress\relax
\let\theasciiabstract\relax\let\theasciikeywords\relax

\let\theasciiemail\relax

\ifplaintex
\font\logobig=cmssbx10 scaled 3836
\font\logomed=cmssbx10 scaled 2557
\else
\font\logobig=cmssbx10 scaled 4200
\font\logomed=cmssbx10 scaled 2800
\fi

\long\def\makeagttitle{   
\count0=\startpage
\agt\hfill      
\hbox to 45truept{\vbox to 0pt{\vglue -13truept{\logomed A\kern -.37em{\logobig 
T}\kern -.38em G}\vss}\hss}
\break
{\small Volume \thevolumenumber\ (\thevolumeyear)
\startpage--\finishpage\nl
Published: \publishdate}

\vglue .25truein

{\parskip=0pt\leftskip 0pt plus
1fil\def\\{\par\smallskip}{\Large\bf\thetitle}\par\medskip} \vglue
0.05truein

{\parskip=0pt\leftskip 0pt plus 1fil\def\\{\par}{\sc\theauthors}
\par\medskip}%
 
\vglue 0.03truein 

{\small\leftskip 25truept\rightskip 25truept{\bf Abstract}\stdspace\theabstract

{\bf AMS Classification}\stdspace\theprimaryclass
\ifx\thesecondaryclass\relax\else; \thesecondaryclass\fi\par
{\bf Keywords}\stdspace \thekeywords\par}\vglue 7truept

}   

\ifplaintex
\font\phead=cmsl9 scaled 950
\font\pnum=cmbx10 scaled 913
\font\pfoot=cmsl9 scaled 950
\headline{\vbox to 0pt{\vskip -4.5mm\line{\small\phead\ifnum
\count0=\startpage ISSN 1472-2739 (on-line) 1472-2747 (printed)
\hfill {\pnum\folio}\else\ifodd\count0\def\\{ }%
\ifx\theshorttitle\relax\thetitle\else\theshorttitle\fi\hfill{\pnum\folio}
\else\def\\{ and }{\pnum\folio}\hfill\ifx\theshortauthors\relax\theauthors
\else\theshortauthors\fi\fi\fi}\vss}}
\footline{\vbox to 0pt{\vglue 0mm\line{\small\pfoot\ifnum\count0=\startpage
\copyright\ \gtp\hfill\else
\agt, Volume \thevolumenumber\ (\thevolumeyear)\hfill\fi}\vss}}
\else
\font=cmsl9 scaled 1050
\font\lnum=cmbx10 
\font=cmsl9 scaled 1050
\makeatletter
\def\@oddhead{{\small\lhead\ifnum\count0=\startpage ISSN 1472-2739 
(on-line) 1472-2747 (printed)\hfill {\lnum\number\count0}\else\ifodd\count0
\def\\{ }\ifx\theshorttitle\relax \thetitle \else\theshorttitle\fi\hfill
{\lnum\number\count0}\else\def\\{ and }{\lnum\number\count0}
\hfill\ifx\theshortauthors\relax 
\theauthors\else\theshortauthors\fi\fi\fi}}\def\@evenhead{\@oddhead}
\def\@oddfoot{\small\lfoot\ifnum\count0=\startpage\copyright\ \gtp\hfill\else
\agt, Volume \thevolumenumber\ (\thevolumeyear)\hfill\fi}
\def\@evenfoot{\@oddfoot}
\makeatother
\fi
\let\maketitlepage\makeagttitle
\let\makeshorttitle\maketitlepage
\let\maketitle\maketitlepage

\newwrite\gtoutfile
\long\gdef\makeheadfile{  
{\def\\{, }\def\s{ }
\immediate\openout\gtoutfile head.xxx
\immediate\write\gtoutfile{Proxy-for: \ifx\theasciiauthors\relax
\theauthors\else\theasciiauthors\fi\s<\ifx\theasciiemail\relax\theemail\else\theasciiemail\fi>}
\immediate\write\gtoutfile{\noexpand\\}
\immediate\write\gtoutfile{Authors: \ifx\theasciiauthors\relax
\theauthors\else\theasciiauthors\fi}
{\def\\{ }\immediate\write\gtoutfile{Title: \ifx\theasciititle\relax
\thetitle\else\theasciititle\fi}}
\immediate\write\gtoutfile{Subj-class: GT or SG, GR etc}
\immediate\write\gtoutfile{MSC-class: \theprimaryclass\ifx\thesecondaryclass\relax\else, \thesecondaryclass\fi}
\immediate\write\gtoutfile{Journal-ref: Algebr. Geom. Topol. \thevolumenumber\s
(\thevolumeyear) \startpage-\finishpage}
\immediate\write\gtoutfile{Comments: Published by Algebraic and
Geometric Topology at}
\immediate\write\gtoutfile{\s\s\s  http://www.maths.warwick.ac.uk/agt/AGTVol\thevolumenumber/agt-\thevolumenumber-\thepapernumber.abs.html}
\immediate\write\gtoutfile{\noexpand\\}
\immediate\write\gtoutfile{}
\ifx\theasciiabstract\relax
\immediate\write\gtoutfile{\theabstract}\else
\immediate\write\gtoutfile{\theasciiabstract}\fi
\immediate\write\gtoutfile{}
\immediate\write\gtoutfile{\noexpand\\}
\immediate\write\gtoutfile{}
\immediate\closeout\gtoutfile}}  

\def\maketitlepage{\makeagttitle\makeheadfile}
\let\makeshorttitle\maketitlepage
\let\maketitle\maketitlepage

\lognumber{37}
\volumenumber{4}
\volumeyear{2004}
\papernumber{37}
\pagenumbers{841}{859}
\received{17 June 2004} 
\accepted{6 October 2004}
\published{7 October 2004}

\usepackage{amsmath,amssymb}
\usepackage{epsf}

\long\def\forget#1\forgotten{}

\newcommand\sg[1]{{\left<{#1}\right>}}
\newcommand\FIGURE[3]{{\begin{figure}[ht!]\cl{\epsfbox{#2}}\caption{#1}\label{#3}\end{figure}}}
\newcommand\FIGUREx[2]{{\begin{figure}[ht!]\cl{\epsfbox{#1}}\nocolon\caption{}\label{#2}\end{figure}}}

\newcommand\set[1]{{\{#1\}}}
\newcommand\eq[1]{{(\ref{#1})}}
\newcommand\eqs[2]{{(\ref{#1})--(\ref{#2})}}
\newcommand\Trip[2]{{{#1}{#2}{#1}={#2}{#1}{#2}}}
\def\limiits{{}} 

\def\C{{\mathbb C}}
\def\CP{{\mathbb{CP}}}

\def\Z{{\mathbb Z}}
\def\P{{\mathbb P}}

\def \Dl{\Delta}

\def \g{\gamma}
\def \G{\Gamma}

\def \s{\sigma}

\def \1{^{-1}}
\def \2{^{-2}}

\def\Pitil{\widetilde\Pi}

\DeclareMathOperator{\Aff}{Aff}
\DeclareMathOperator{\Gal}{Gal}
\DeclareMathOperator{\Ker}{Ker}

\newtheorem{theorem}{Theorem}[section] 
\newtheorem{tab}[theorem]{Table}

\newtheorem{proposition}[theorem]{Proposition}
\newtheorem{corollary}[theorem]{Corollary}
\newtheorem{lemma}[theorem]{Lemma}

\newtheorem{defn}[theorem]{Definition}

\begin{document}

\covertitle{Higher degree Galois covers of ${\noexpand\bf CP}^1\times T$}

\asciititle{Higher degree Galois covers of CP^1 x T}

\title[Galois covers of $\C\P^1\times T$]{Higher degree Galois covers of 
$\C\P^1\times T$}

\authors{Meirav Amram\\David Goldberg}
\addresses{Einstein Institute for Mathematics, the Hebrew University, 
Jerusalem, Israel\\\smallskip\\
Mathematics Department, Colorado State University, Fort Collins, CO
80523 USA}
\asciiaddress{Einstein Institute for Mathematics, the Hebrew University, 
Jerusalem, Israel\\
Mathematics Department, Colorado State University, Fort Collins, CO
80523 USA}

\asciiemail{ameirav@math.huji.ac.il, david_j_goldberg@hotmail.com}
\gtemail{\mailto{ameirav@math.huji.ac.il}{\qua\rm and\qua}\href{mailto:david_j_goldberg@hotmail.com}{david\_j\_goldberg@hotmail.com}}

\asciiabstract{%
Let T be a complex torus, and X the surface CP^1 x T.  If T is
embedded in CP^{n-1} then X may be embedded in CP^{2n-1}.  Let X_Gal
be its Galois cover with respect to a generic projection to CP^2. In
this paper we compute the fundamental group of X_Gal, using the
degeneration and regeneration techniques, the Moishezon-Teicher braid
monodromy algorithm and group calculations. We show that pi_1(X_Gal) =
Z^{4n-2}.}

\begin{abstract}
Let $T$ be a complex torus, and $X$ the surface $\C\P^1 \times T$.
If $T$ is embedded in $\C\P^{n-1}$ then $X$ may be embedded in
$\C\P^{2n-1}$.  Let $X_{\Gal}$ be its Galois cover with respect
to a generic projection to $\C\P^2$. In this paper we compute the
fundamental group of $X_{\Gal}$, using the degeneration and
regeneration techniques, the Moishezon--Teicher braid monodromy
algorithm and group calculations. We show that $\pi_1(X_{\Gal}) = \Z^{4n-2}$.
\end{abstract}

\asciikeywords{Galois cover, fundamental group, 
generic projection, Sieberg-Witten invariants}

\maketitle

\section{Overview}\label{overview}

Let $T$ be a complex torus embedded in $\C\P^{n-1}$.  The surface
$X = \C\P^1\times T$ can be embedded into projective space using the
Segre embedding from $\C\P^1\times\C\P^{n-1}\to\C\P^{2n-1}$.
We compute the fundamental group of the Galois cover of $X$ with respect
to a generic projection $f$ from $X\subset\C\P^{2n-1}$ to $\C\P^2$.
This map has degree $2n$.  The Galois cover can be defined as the closure
of the $2n$--fold fibered product
$X_{\Gal}=\overline{X\times_f\cdots\times_fX-\Delta}$ where
$\Delta$ is the generalized
diagonal.  The closure is necessary because the branched fibers
are excluded when $\Delta$ is omitted.

Since the induced map $X_{\Gal}\to\C\P^2$ has the same branch
curve $S$ as $f\co X\to\C\P^2$, the fundamental group
$\pi_1(X_{\Gal})$ is related to $\pi_1(\C\P^2-S)$. In fact it is a
normal subgroup of $\Pitil_1$, the quotient of $\pi_1(\C\P^2-S)$ by the normal
subgroup generated by the squares of the standard generators. In
this paper we employ braid monodromy techniques, the van Kampen theorem
and various computational methods of groups to compute a
presentation for the quotient $\Pitil_1$ from which
$\pi_1(X_{\Gal})$ can be derived. Our main result is that
$\pi_1(X_{\Gal}) = \Z^{4n-2}$ (Theorem \ref{main}).  This extends a
previous result for $n=3$ proven in \cite{AGTV}.

The paper follows the structure of \cite{AGTV}. In Section \ref{degen} we
describe the degeneration of the surface $X=\C\P^1\times T$ and
the degenerated branch curve. In Section \ref{regen}
we regenerate the branch curve and its braid
monodromy factorization to get a presentation for
$\pi_1(\C^2-S)$, the fundamental group of the complement of
the regenerated branch curve in $\C^2$.
In Section \ref{XGal} we compute $\pi_1(X_{\Gal})$ as the kernel of a
permutation monodromy map using the Reidmeister--Schreier method.

\section{Degeneration of  $\C\P^1\times T$}\label{degen}

To compute the braid monodromy of the branch curve $S$, we degenerate $X$ to
a union of projective planes $X_0$.  The branch curve degenerates to a
union of lines $S_0$ which are the images of the intersections of the
planes of $X_0$.  We use the following degeneration:
Embedded as a degree $n$ elliptic normal curve in $\C\P^{n-1}$, the
torus $T$ degenerates to a cycle of $n$~projective lines.  Under the
Segre embedding $X$ degenerates to a cycle of $n$~quadrics,
$Q_i \cong \CP^1\times \CP^1$.  Call this space $X_1$, see Figure
\ref{nSQ}.

\FIGURE{The space $X_1$}{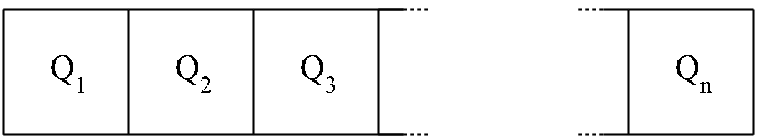}{nSQ} 

Each quadric in $X_1$ shares a projective line with its two neighbors.
As a cycle of quadrics, it is understood that $Q_1$ and $Q_n$ intersect
as well, though it is not clearly indicated in Figure \ref{nSQ}.  Hence
the left and right edges should be identified to make an $n$--prism.
Each quadric in $X_1$ can be further degenerated to a union of two
projective planes.  In Figure \ref{scomp} this is represented by a
diagonal line which divides each square into two triangles, each
representing $\CP^2$.
We shall refer to this diagram as the simplicial complex of $X_0$.

\FIGURE{The simplicial complex $X_0$}{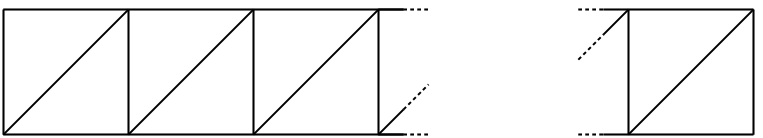}{scomp}

A common edge between two triangles represents the intersection
line of the two corresponding planes.  The union of the
intersection lines in $X_0$ is the ramification curve of
$f_0\co X_0\to \C\P^2$, denoted by $R_0$. Let  $S_0 = f_0(R_0)$
be the degenerated branch curve.  It is a line arrangement,
composed of the images of all the intersection lines.

Each vertex of the simplicial complex represents an intersection
point of three planes.  These are the singular points of
$R_0$. Each of these vertices is called a $3$--point (reflecting the
number of planes which meet there).

The vertices may be given any convenient enumeration.  We have chosen
left to right, bottom to top enumeration, see Figure \ref{vnum}.
Because the left and right edges are identified, so are the
corresponding vertices.

\FIGURE{The enumeration of vertices}{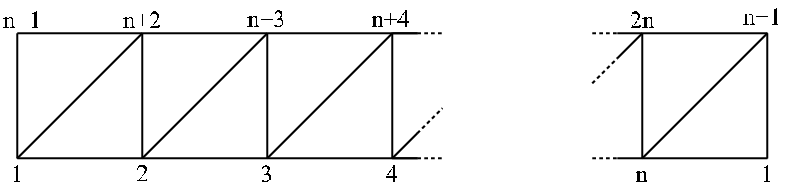}{vnum}

In order to use a result concerning the monodromy of dual to generic
line arrangements \cite{Mo} we number the edges based upon the enumeration
of the vertices using reverse lexicographic ordering: if $L_1$ and
$L_2$ are two lines with end points $\alpha _1, \beta _1$ and
$\alpha _2, \beta _2$ respectively $(\alpha _1 < \beta _1, \alpha
_2 < \beta _2)$, then $L_1 < L_2$ iff $\beta_1 < \beta _2$, or
$\beta _1 = \beta _2$ and $\alpha _1 < \alpha _2$. The resulting
enumeration is shown in Figure \ref{enum}.  The horizontal lines at the
top and bottom do not represent intersections of planes and hence are not
numbered.

\FIGURE{The enumeration of lines}{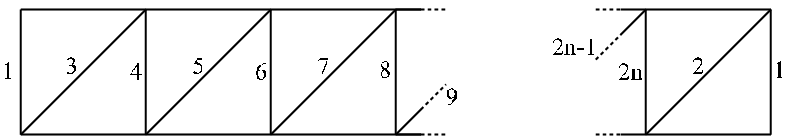}{enum}

The triangles can be numbered in any order, so we choose an enumeration
which will simplify future computations.

\FIGURE{The enumeration of planes}{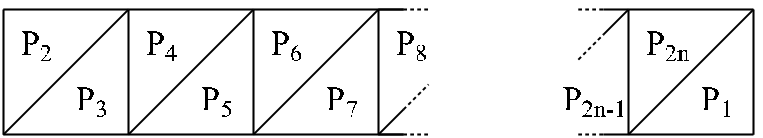}{pnum}

The braid monodromy of $S_0$ is easily computed since it is a subset of
a dual to generic line arrangement \cite{Mo}.  From this, the braid
monodromy of the full branch curve $S$ can be regenerated.

\section{Regeneration of the branch curve}\label{regen}

The degeneration of $X$ through $X_1$ to $X_0$ takes place in $\C\P^{2n-1}$
so at every step of the process the generic projection
$\C\P^{2n-1}\to\C\P^2$ restricts to a map $f_i\co X_i\to\C\P^2$.
Each map $f_i$ has its own branch curve $S_i$.
Starting from the degenerated branch curve $S_0$, we reverse the steps of
the degeneration of $X$ to regenerate the braid monodromy of $S$.

\subsection{The braid monodromy of $S_0$}\label{BMS0}

We have enumerated the $2n$ planes $P_1,\ldots,P_{2n}$ which comprise
$X_0$, the $2n$ intersection lines $\hat L_1,\ldots,\hat L_{2n}$ which
comprise $R_0$, and their $2n$ intersection points
$\hat V_1,\ldots,\hat V_{2n}$.  Let
$L_i$ and $V_k$ denote the projections of $\hat L_i$ and $\hat
V_k$ to $\C\P^2$ by the map $f_0$.  Clearly the degenerated branch curve
$S_0=\bigcup \limiits_{i=1}^{2n}L_i$.  Each of the $V_k$ is an
intersection point of $S_0$.  Additionally, every pair of lines
$\hat L_i$ and $\hat L_j$ which do
not intersect in $R_0$ must have a simple intersection when projected to
$S_0\subset\C\P^2$.
Thus the braid monodromy of $S_0$ consists of cycles $\Dl^2_k$
corresponding to the $V_k$, and full twists $D_{ij}$ corresponding to
the simple intersection $L_i\cap L_j$. 

Since $S_0$ is a sub-arrangement of a dual to generic arrangement the
precise forms of the braids $\Dl^2_k$ and $D_{ij}$ can be found using
Moishezon's result \cite{Mo}, summarized in Theorem IX.2.1 of
\cite{MoTe2}.

The braids take place in a generic fiber $\C_u$ of the projection
$\pi\co\C\P^2\to\C\P^1$.  For each line we will refer to the intersection
$L_i\cap\C_u$ as the point~$i$.  These are the points which are braided
in $\C_u$.  Since they are so closely related we
may often use the concepts of the line $L_i$ and the point $i$
interchangeably.

Recall also that each line $L_i$ is initially denoted by a pair of numbers
indicating which two vertices it contains.  The lines are sorted based
on the second vertex then the first, producing the single index $i$.
We will need to refer to the pairs of vertices associated to the lines 
$L_i$ and $L_j$ in defining the braids $D_{ij}$.

Assume $i<j$.  From \cite{Mo}, \cite{MoTe2} we know that each $D_{ij}$ is a
particular full twist of the points
$i$ and $j$ in which $i$ is brought next
to $j$ by passing over most of the intervening points, but under those
points which share the same second vertex as $j$.  The path described is
denoted by $\tilde z_{ij}$ and the corresponding half-twist is called
$\tilde Z_{ij}$ so we may say that $D_{ij}=\tilde Z_{ij}^2$.  In the context
of our enumeration of line this means that when $j$ is a vertical
(see Figure \ref{enum}) then the
path passes under the preceding diagonal.  With the exception of lines $1$
and $2$, the odds are diagonal and the evens are vertical.
Representative examples are shown in Figure \ref{Dij}.

\FIGURE{Full twists $D_{ij}$}{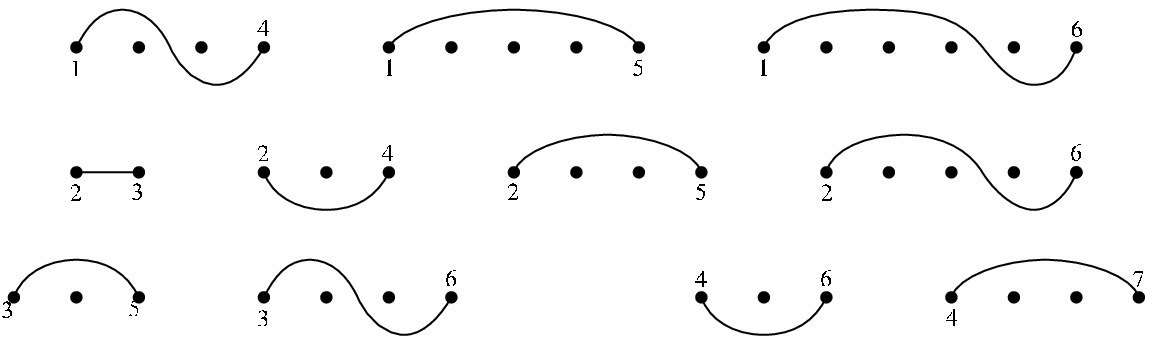}{Dij}

Each $\Dl^2_k$ is a particular full twist of all the points corresponding to
lines through $V_k$ in which the points pass under any intervening points.
In our example each $V_k$ is the intersection of exactly $2$~lines.
Let $i$ and $j$ be the indices of the two lines through $V_k$.  Let
$\underline z_{ij}$ denote the path from $i$ to $j$ under all the points
in between and $\underline Z_{ij}$ the corresponding half-twist.  In
this terminology $\Dl^2_k=\underline Z_{ij}^2$.
Since most of the pairs of lines meeting at the $V_k$
are consecutively numbered there aren't many intervening points to worry
about.  The exceptions are lines $1$ and $3$ intersecting at $V_1$ and
lines $2$ and $2n$ intersecting at $V_n$.

Using the braid monodromy of $S_0$ as a template, the braid monodromy of
$S$ can be obtained according to a few simple regeneration rules.

\subsection{The braid monodromy of $S$}\label{BMS}

When $X_0$ is regenerated to $X$, each of the lines $L_i$ in $S_0$
divides into two sheets of the branch curve $S$.  So in the generic
fiber $\C_u$ the point $i$ divides into two points which we shall call
$i$ and $i'$.  As each intersection point of $S_0$ splits into a
collection of singularities of $S$, the associated braids $D_{ij}$ and
$\Dl_k^2$ also split into collections of braids in predictable ways.
(Basic regeneration rules are proven in \cite{MoTe4}.  Application to
the specific types of singularities found here is as in \cite{AGTV}.)

Assume $i<j$.  For each pair of lines $\hat L_i$ and $\hat L_j$ which do
not intersect in $R_0$ there is a simple node in $S_0$ with monodromy
$D_{ij}$ where $L_i$ and $L_j$ intersect.  When $L_i$ and $L_j$ divide
into two sheets the resulting figure has four nodes, each with its own
monodromy.  So each $D_{ij}=\tilde Z^2_{ij}$ becomes the four braids:
$\tilde Z^2_{ij}$, $\tilde Z^2_{ij'}$, $\tilde Z^2_{i'j}$, and
$\tilde Z^2_{i'j'}$, which we summarize with the symbol
$\tilde Z^2_{ii',jj'}$.  Figure \ref{Ziijj} shows two representative
examples.

\FIGURE{Regenerated collections of braids $\tilde Z^2_{11',44'}$ and $\tilde Z^2_{22',55'}$}{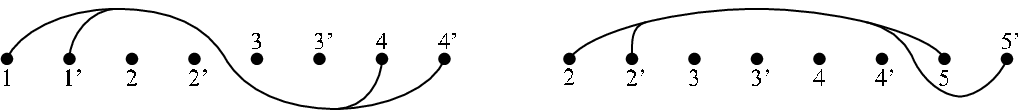}{Ziijj}

Each $\hat V_k$ is the intersection of three planes in $X_0$.  For that
reason this type of singularity is called a $3$--point.  These points
regenerate in two steps.  First in $X_1$ pairs of planes are regenerated
to quadrics so in the branch curve $S_1$ each diagonal line becomes a conic,
tangent to the adjacent vertical lines.  Near each $3$--point this
creates a branch point and a point of tangency.  Then when $X$ is
regenerated, the vertical lines divides, and the points of tangency
become three cusps each.  Hence the braid $\Dl^2_k$ becomes braids.
The branch point yields a half-twist, and each cusp yields a $3/2$-twist.
Let $i$ and $j$ respectively be the indices of the vertical and diagonal
lines meeting at $V_k$.  The symbol $Z_{jj'(i)}$ denotes the half-twist
and $Z^3_{ii',j'}$ denotes the three $3/2$-twists.  An illustrative
example is provided in Figure \ref{threept}

\FIGURE{Regenerated braids for $V_1$: $Z^3_{11',3'}$ and $Z_{33'(1)}$}{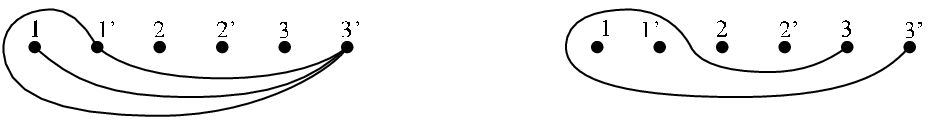}{threept}

These braids comprise the full braid monodromy of the branch curve $S$.
To check this recall that the product, suitably ordered, of all the
monodromy gives a factorization of $\Dl^2_S$.  Since $S$ has degree $4n$,
$\Dl^2_S$
has degree $4n(4n-1) = 16n^2-4n$.  In the simplicial complex there are
$2n$ different $3$--points.  Each produces $3$~cusps and a branch point
for a combined degree of $2n\cdot(3\cdot 3 + 1) = 20n$.  There are also
$2n^2-3n$ incidental intersections which do not appear in the simplicial
complex.  Each produces $4$~nodes for a combined degree of
$(2n^2-3n)(4\cdot 2) = 16n^2 - 24n$.  So together they have the
required total degree of $16n^2-4n$.

\subsection{The fundamental group of $\C^2-S$}\label{pi1c}

Let $X^{\Aff}$ denote $f^{-1}(\C^2)$, the portion of $X$ lying above some
affine $\C^2\subset\C\P^2$.  We will use the van Kampen Theorem \cite{vK}
to produce a presentation for the fundamental group of the affine complement
of the branch curve $\pi_1(\C^2-S)$ from the braid monodromy
factorization computed in Subsection \ref{BMS}.  Since the braid
monodromy factorization is only defined up to Hurwitz equivalence of
factorizations it is clear that equivalent factorization must give the
same fundamental group.  So any relations obtained from equivalent
factorizations may be included without changing the group.  To get extra
relations
we will be making liberal use of the following two invariance theorems.

\begin{theorem}[Invariance Theorem]\label{thi}
The braid monodromy factorization of $S$
is invariant under conjugation by any product of half-twists of the form
$Z_{jj'}$.
\end{theorem}

Note that the braids $Z_{ii'}$ and $Z_{jj'}$ commute for all $i$ and $j$
because the path from $i$ to $i'$ doesn't intersect the path from $j$ to $j'$.
Theorem \ref{thi} is proven in \cite{MoTe4} and \cite{AGTV}.  Philosophically
it is true because going around a neighborhood of the degenerated surface
$X_0$ causes the two sheets associated to $L_j$ to interchange, which
conjugates all braids by $Z_{jj'}$.  The $Z_{jj'}$ can be used
independently because in this figure the lines $L_i$ can be regenerated
independently.

\begin{theorem}[Conjugation Theorem]\label{thc}
The braid monodromy factorization of $S$
is invariant under complex conjugation in $\C_u$ if the order of factors
is reversed.
\end{theorem}

This theorem is proven in \cite{MoTe4}.  Philosophically it is true
because complex conjugation in two dimensions is orientation preserving.
Reversing the order of factors corresponds to complex conjugation in the
other dimension.

By the van Kampen Theorem \cite{vK}, there is a surjection from the
fiber $\pi_1(\C_u-\{j,j'\})$ onto  $\Pi_1=\pi_1(\C^2-S)$.
The fundamental group of the fiber is
freely generated by $\{\G_j,\G_{j'}\}_{j=1}^{2n}$, where $\G_j$ and
$\G_{j'}$ are loops in $\C_u$ around $j$ and $j'$ respectively.
These loops are explicitly constructed in \cite{AGTV}.
Thus the images of $\{\G_j,\G_{j'}\}_{j=1}^{2n}$ are generators for $\Pi_1$. 
Without too much confusion we will refer to the images as
$\{\G_j,\G_{j'}\}_{j=1}^{2n}$ as well.  Each braid in the braid monodromy
factorization of $S$ induces a relation on $\Pi_1$ through its
natural action on $\C_u-\{j,j'\}_{j=1}^{2n}$ \cite{vK}.

Using the techniques described above we get a presentation for
the affine complement of the branch curve $S$.  For simplicity of
notation we will be using the following shorthand:

$\G_{(i)}$ stands for all of the conjugates of $\G_i$ by integer powers of
$Z_{ii'}$.  These include
$...\G_{i'}^{-1}\G_i^{-1}\G_{i'}\G_i\G_{i'},\ \ 
          \G_{i'}^{-1}\G_i\G_{i'},\ \ \G_{i'},\ \ \G_{i},\ \ 
          \G_i\G_{i'}\G_i^{-1},\ \ 
          \G_i\G_{i'}\G_i\G_{i'}^{-1}\G_i^{-1}...$

$\G_{ii'}$ stands for either $\G_i$ or $\G_{i'}$.

$\G_{\check ii'}$ stands for either $\G_{i'}^{-1}\G_i\G_{i'}$ or $\G_{i'}$.

$\G_{j\hat j}$ stands for either $\G_j$ or $\G_j\G_{j'}\G_j^{-1}$.

\begin{theorem}\label{pres1}
The group $\Pi_1$ is generated by $\set{\G_j, \G_{j'}}_{j=1}^{2n}$ with
the following relations:
\begin{eqnarray}
{}[\G_{(i)},\G_{(j)}] & = & 1 \qquad \mbox{if the lines $i,j$ are
disjoint in $X_0$} \label{commeq1}\\
{}\G_{(i)}\G_{(j)}\G_{(i)}  & = & \G_{(j)}\G_{(i)}\G_{(j)} \qquad
\mbox{if the lines $i,j$ intersect}\label{tripeq1}\\
{}\G_{j'}^{-1}\G_i\G_{i'}\G_j\G_{i'}^{-1}\G_i^{-1} & = & 1 \qquad
\mbox{if vertical $i$ intersects diagonal $j$.} \label{bpeq1}
\end{eqnarray}
\end{theorem}

\begin{proof}
The relations \eq{bpeq1} come from the braids $Z_{jj'(i)}$.  The
relations from $Z^3_{ii',j'}$ include
$\G_{i'}\G_{j'}\G_{i'}=\G_{j'}\G_{i'}\G_{j'}$ and two other variants.
Using the Invariance Theorem \ref{thi} we get the triple relations
\eq{tripeq1} in their full generality.  It remains to prove commutation
relations \eq{commeq1}.

For $i<j$, if $j$ is odd then $j$ is diagonal so $\tilde Z^2_{ij}$ only goes
over the intervening points.  As a result the complex conjugates of the
regenerated braids give relations $[\G_{ii'},\G_{j\hat j}]=1$, see 
Figure \ref{overcomm}.  These four relations can easily be seen to
generate the rest of $[\G_{(i)},\G_{(j)}]=1$.

\FIGUREx{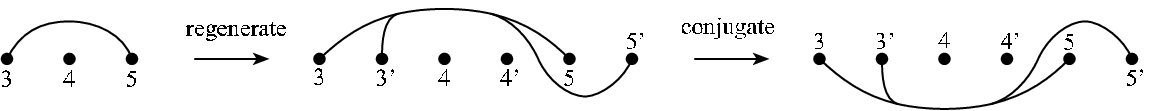}{overcomm}

If $j$ is even and $i=j-2$ then $\tilde Z^2_{ij}$ only goes under the
intervening points.  As a result the regenerated braids give relations
$[\G_{\check ii'},\G_{jj'}]=1$, see Figure \ref{undercomm}.  Again these
four relations generate all of $[\G_{(i)},\G_{(j)}]=1$.

\FIGUREx{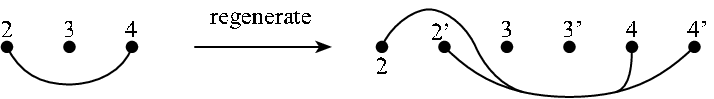}{undercomm}

For the remainder of the commutators when $j$ is even $\tilde Z^2_{ij}$ goes
under $j-1$ and over all the other intervening points.  The complex
conjugates of the regenerated braids give relations
$[\G_{ii'},\G_{j-1}\G_{j-1'}\G_{j\hat j}\G_{j-1'}^{-1}\G_{j-1}^{-1}]=1$, 
see Figure \ref{conjcomm}.  Since $j-1$ is odd it normally follows from the
arguments above that all $\G_{(i)}$ commute with all $\G_{(j-1)}$.  Hence
the commutator above simplifies to $[\G_{ii'},\G_{j\hat j}]=1$ which
implies $[\G_{(i)},\G_{(j)}]=1$.

     The one exception is when $i=1$ and $j=4$ because $\G_{(1)}$ and
$\G_{(3)}$ satisfy triple relations rather than commutators.  In this
case, the regenerated braids from $\tilde Z^2_{14}$ without complex
conjugation give
$[\G_{2'}^{-1}\G_2^{-1}\G_{\check 11'}\G_2\G_{2'},\G_{44'}]=1$.  Now the
fact that $\G_{(2)}$ and $\G_{(4)}$ commute will finish the proof.
\end{proof}

\FIGUREx{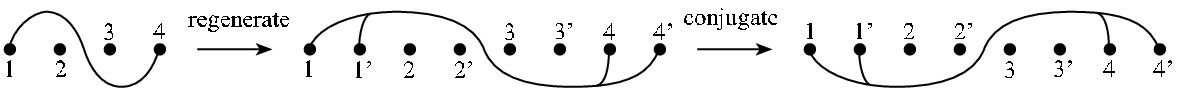}{conjcomm}

\section{The Galois cover of $X$}\label{XGal}

\subsection{The homomorphism $\psi$}\label{hompsi}

$X^{\Aff}-f^{-1}(S)$ is a degree $2n$ covering space of $\C^2-S$.
Let
$$\psi\co \pi_1(\C^2-S)\to S_{2n}$$
be the permutation monodromy of this cover.  We can compute $\psi$
precisely by considering the degeneration of $X$ to $X_0$.  In $X_0$ the
sheets are the planes $P_k$ which we numbered in Figure \ref{pnum}.  In
$X_0$ there is no monodromy since $X_0 - R_0$ breaks up into a disjoint
union so it is impossible to get from one sheet to another.

For a regeneration of $X$ near $X_0$ the sheets of $X$ are very close to
the planes of $X_0$, so we can use the same numbering.  As we have seen
already, the ramification curve $R$ locally has two pieces near each
line $\hat L_i$ of $R_0$.  The regeneration of the monodromy here looks
very much like a node pulling apart into two branch points.  In fact
when restricted to the generic fiber $\C_u$ it is precisely the node $i$
pulling apart into the two branch points $i$ and $i'$ with identical
monodromy.  If $\hat L_i$ is the intersection of $P_k$ and $P_\ell$ then
it is clear that $\psi(\G_i)=\psi(\G_{i'})=(k~\ell)$.
Based on the enumerations of lines and planes in
Figures \ref{enum} and \ref{pnum} we have:
\begin{defn}\label{defpsi}
The map $\psi\co\pi_1(\C^2-S)\to S_{2n}$ is given by
\begin{eqnarray*}
\psi(\G_1)=\psi(\G_{1'})&=&(1~2),\\
\psi(\G_i)=\psi(\G_{i'})&=&(i{-}1~i)\qquad \mbox{for $3\le i\le 2n$},\\
\psi(\G_2)=\psi(\G_{2'})&=&(2n~1)
\end{eqnarray*}
\end{defn}

The reader may wish to check that $\psi$ is well defined by testing
the relations given in Theorem \ref{pres1}), but this is of course
guaranteed by the theory. From the definitions \ref{defpsi} we see that
$\psi$ is surjective.  In fact the images of $\G_2,\cdots,\G_{2n}$
generate $S_{2n}$.

\subsection{The fundamental groups of $X^{\Aff}_{\Gal}$ and $X_{\Gal}$ in terms of $\psi$}\label{twogroups}

The sheets of the Galois cover of $X$ are labeled by permutations of
the $2n$ sheets of $X$.  Therefore an element of $\pi_1(\C^2-S)$
lifts to a closed path in $X^{\Aff}_{\Gal}$ if and only if it has no
permutation monodromy, which means it is in the kernel of $\psi$.
Among these are elements $\G_j^2$ and $\G_{j'}^2$ which lift to double
loops around simple ramification curves.  These double loops vanish in
$\pi_1(X^{\Aff}_{\Gal})$.  We may define a quotient
\begin{equation}\label{Pitildef}
\Pitil_1=\frac{\pi_1(\C^2-S)}{\sg{\G_j^2,\G_{j'}^2}}
\end{equation}
in which the elements $\G_j^2$ and $\G_{j'}^2$ are forcibly killed.
Since these elements are already in the kernel, $\psi$ remains well defined on
$\Pitil_1$.  Let $\mathcal{A}$ be the kernel of
$\psi\co \Pitil_1\rightarrow S_{2n}$.  We have a short exact sequence
sequence
\begin{equation}\label{short}
1 \longrightarrow \mathcal{A}
\longrightarrow \Pitil_1\stackrel{\psi}{\longrightarrow} S_{2n}\rightarrow 1.
\end{equation}

With these observations it is clear that
the fundamental group $\pi_1(X^{\Aff}_{\Gal})$ is isomorphic to
$\mathcal{A}$.  Based on the presentation of $\Pi_1$ given in Theorem
\ref{pres1} we can write a presentation for $\Pitil_1$

\begin{theorem}\label{pres2}
The group $\Pitil_1$ is generated by $\set{\G_j, \G_{j'}}_{j=1}^{2n}$ with
the following relations:
\begin{eqnarray}
{}\G_i^2 & = & 1 \label{sqeq}\\
{}\G_{i'}^2 & = & 1 \label{sqeq2}\\
{}[\G_{(i)},\G_{(j)}] & = & 1 \qquad \mbox{if the lines $i,j$ are
disjoint in $X_0$} \label{commeq2}\\
{}\G_{(i)}\G_{(j)}\G_{(i)}  & = & \G_{(j)}\G_{(i)}\G_{(j)} \qquad
\mbox{if the lines $i,j$ intersect}\label{tripeq2}\\
{}\G_{j'}\G_i\G_{i'}\G_j\G_{i'}\G_i & = & 1 \qquad
\mbox{if vertical $i$ intersects diagonal $j$}. \label{bpeq2}
\end{eqnarray}
$\Gamma_{(i)}$ stands for any odd length word in the
infinite dihedral group $\sg{\Gamma_{i},\Gamma_{i'}}$.
\end{theorem}

    This result can be extended to the whole surface $X_{\Gal}$ by
including the projective relation.  $\pi_1(\C\P^2-S)$ has the same
generators and relations as the affine $\pi_1(\C^2-S)$ with one
additional relation:
$$\G_1\G_{1'}\G_2\G_{2'}\cdots\G_{2n}\G_{2n'} = 1.$$
This relation is in the kernel of $\psi$ so it does not interfere with
the definition of $\psi$.  If we add the projective relation to the
presentation of $\Pitil_1$ then we can compute $\pi_1(X_{\Gal})$ as the
kernel of $\psi$ as before.

\subsection{The splitting $\varphi$ and the Reidmeister--Schreier method}\label{splitRMS}

It will be useful to have a splitting of the short exact sequence
\eq{short}.  Define the splitting map $\varphi\co S_{2n}\to\Pitil_1$
as follows:
\begin{defn}\label{defphi}
The map $\varphi\co S_{2n}\to\Pitil_1$ is given by
\begin{eqnarray*}
\varphi(2n~1)    &=& \G_2,\\
\varphi(i{-}1~i) &=& \G_i \qquad \mbox{for $3\le i\le 2n$}.
\end{eqnarray*}
\end{defn}
From the definitions of $\psi$ and $\varphi$ it is clear that
$\psi\circ\varphi$ is the identity on $S_{2n}$, so it remains to check
that $\varphi$ is well defined.  The set of generators
$(2~3),\ (3~4),\ \cdots,\ (2n{-}1~2n),\ (2n~1)$ used in the definition
of $\varphi$ satisfy the following relations:  The square of each
generator is the identity;  Non intersecting generator commute; 
Intersecting generators satisfy triple relations.  These relations 
are suffice (it is well known, as it was proven in \cite{AGTV}). 
Note that all of these relations are also satisfied by $\G_3,\ 
\G_4,\ \cdots,\ \G_{2n},\ \G_2$ in $\Pitil_1$.  
In fact they are in the presentation as relations \eq{sqeq},
\eq{commeq2}, and \eq{tripeq2}.

We use the Reidmeister--Schreier method to find a presentation for
the kernel $\mathcal{A}$ of the map $\psi\co\Pitil_1
\rightarrow S_{2n}$.  Since $\psi$ is split by $\varphi$ we can write the
generators of $\mathcal{A}$ as follows:
$$\g(\s,\G)=\s\G(\varphi\psi(\G))^{-1}\s^{-1}\quad
\forall\s\in\varphi(S_{2n}),\quad
\mbox{$\forall\G$ generator of $\Pitil_1$}.$$
To simplify notation let $\overline\G=\varphi\psi(\G)$ denote the
projection of $\G$ onto the image of $\varphi$.
Using this notation we get generators:
$$\g(\s,\G)=\s\G\overline\G^{-1}\s^{-1}.$$
Since $S_{2n}\cong\varphi(S_{2n})\subset\Pitil_1$ we will not
distinguish between the two groups.  We will think of $S_{2n}$ as a
subgroup of $\Pitil_1$, so $\s\in S_{2n}$ above are permutations.

The relations of $\Pitil_1$ can be translated into expressions in these 
generators by the following process.
If the word
$\omega = \G_{i_1}\G_{i_2}\cdots\G_{i_t}$ represents an element of $\Ker\psi$
then $\omega$ can be rewritten as the product
$$\tau (\omega) = \g(1,\G_{i_1}) \g(\overline\G_{i_1},\G_{i_2})
\cdots
\g(\overline{\G_{i_1} \cdots \G_{i_{t-1}}}, \G_{i_t}).$$

\begin{theorem}[Reidmeister--Schreier]\label{RMS}
Let $\set{R}$ be a complete set of relations for $\Pitil_1$.
Then $\mathcal{A} = \Ker \psi$ is generated by the $\g(\s,\G)$ with
the relations $\{\tau(\s r\s^{-1})\}_{r\in R, \s\in S_{2n}}$.
\end{theorem}

We use this method to find generators and relations for the kernel.

\subsection{Generators for $\mathcal{A} = \Ker \psi$}\label{genker}

By Theorem \ref{RMS}, $\mathcal{A}$ is generated
by the elements $\s\G_j\overline\G_j^{-1}\s^{-1}$ and
$\s \G_{j'}\overline\G_{j'}^{-1}\s^{-1}$, $1 \leq j \leq 2n$, $\s\in S_{2n}$.
We compute $\overline\G_j$ and $\overline\G_{j'}$.
Recall that $\sg{\G_2,\dots,\G_{2n}} = S_{2n}$ is the
image of $\varphi$, so for $j \neq 1$
we get $\overline\G_j = \overline\G_{j'} = \G_j$, and the
associated generators are
\begin{equation}\label{defAsj}
A_{\s,j}=\s\G_{j}\G_{j'}\s^{-1}.
\end{equation}

The permutation $(1\,2)$ can be expressed in terms of the generators
of $S_{2n}$ corresponding to $\G_2,\dots,\G_{2n}$ as follows:
$$(1\,2) = (2n\,1)\cdots(4\,5)(3\,4)(2\,3)(3\,4)(4\,5)\cdots(2n\,1).$$
So for $j = 1$ we have that
$\overline\G_1=\overline\G_{1'}=
   \G_2\G_{2n}\cdots\G_4\G_3\G_4\cdots\G_{2n}\G_2$.

Since we are identifying $S_{2n}$ with
$\varphi(S_{2n})=\sg{\G_2, \ldots, \G_{2n}}$ it is reasonable to write the
above as simply
$$\overline\G_1=\overline\G_{1'}=(1\,2).$$
So we get generators
\begin{eqnarray}
X_{\s} &=& \s(1\,2)\G_1\s^{-1}\label{Xdef}\\
B_{\s} &=& \s(1\,2)\G_{1'}\s^{-1}.\label{Bdef}
\end{eqnarray}

Since $X_\s^{-1}B_\s =
\s \G_1 \G_{1'} \s^{-1}$ we can define $A_{\s,1}=X_\s^{-1}B_\s$
for $j=1$ to get the following result:
\begin{corollary}
The group $\mathcal{A} = \Ker \psi$ is generated by
$A_{\s,j}$, $X_\s$, for $\s\in S_{2n}$ and $j = 1,\dots,2n$.
\end{corollary}

\subsection{Reducing the set of generators for $\mathcal{A}$}\label{reduce}

First we show that the $A_{\s,j}$ are not needed for $j = 2,\dots,2n$.

\begin{theorem}\label{As1Xs}
$\mathcal{A}$ is generated by $\{A_{\s,1}, X_\s\}$.
\end{theorem}
\begin{proof}
The following relations are translations of the relations \eq{bpeq2} from
Theorem \ref{pres2}.  Together they show that all of the
$A_{\s,j}$ for $j\ne 1$ can be written in terms of the $A_{\s,1}$.
Derivations follow the table.

\begin{tab}\label{Vtab-rel}
Translations of the branch point relations:
\begin{eqnarray}
A_{\s,3}  & = & A_{\s(23),1}A_{\s,1}^{-1} \label{AA1}\\
A_{\s,3}  & = & A_{\s(234),4} \label{Aconj3}\\
A_{\s,4}  & = & A_{\s(345),5}\\
\null   &\vdots& \null \notag \\
A_{\s,2n-1} &=& A_{\s(2n-2\,2n-1\,2n),2n}\\
A_{\s,2n} & = & A_{\s(2n-1\,2n\,1),2}\\
A_{\s,2}  & = & A_{\s(2n\,1),1}A_{\s,1}^{-1}.\label{AA2}
\end{eqnarray}
\end{tab}

We use the relations \eq{bpeq2} of Theorem \ref{pres2}. Let $I$ denote the
identity element of $S_{2n}$, so that by definition
$A_{I,j} = \Gamma_j\Gamma_{j'}$.

To prove \eq{AA1} we use $V_1$ which has diagonal $j=3$ and vertical $i=1$.
The branch point relation is
$1 = \G_{3'}\G_1\G_{1'}\G_3\G_{1'}\G_1=
(\G_{3'}\G_3)(\G_3\G_1\G_{1'}\G_3)(\G_{1'}\G_1)=
A_{I,3}^{-1}A_{(23),1}A_{I,1}^{-1}$ so we get
$A_{I,3}=A_{(23),1}A_{I,1}^{-1}$.

To prove \eq{AA2} we use $V_{n+1}$ which has diagonal $j=2$ and vertical
$i=1$.  The relation is
$1=\G_{2'}\G_1\G_{1'}\G_2\G_{1'}\G_1=
(\G_{2'}\G_2)(\G_2\G_1\G_{1'}\G_2)(\G_{1'}\G_1)=
A_{I,2}^{-1}A_{(2n\,1),1}A_{I,1}^{-1}$; so we get
$A_{I,2}=A_{(2n\,1),1}A_{I,1}^{-1}$.

To prove \eq{Aconj3} we use $V_{n+2}$ which has diagonal $j{=}3$ and
vertical $i{=}4$.  The corresponding relation is
$1=\G_{3'}\G_4\G_{4'}\G_3\G_{4'}\G_4=
(\G_{3'}\G_3) \G_3\G_4\G_{4'}\G_3\G_{4'}\G_4 =
(\G_{3'}\G_3) \G_3\G_4\G_{3}\G_{4'}\G_{3}\G_4 =
(\G_{3'}\G_3) \G_4\G_3\G_{4}\G_{4'}\G_{3}\G_4 =
A_{I,3}^{-1} A_{(34)(23),4}$.  As a result $A_{I,3}=A_{(234),4}$.

The rest of Table \ref{Vtab-rel} follows in much the same manner as
\eq{Aconj3}.
\end{proof}

Combining all of the relations in Table \ref{Vtab-rel}, we obtain one 
relation among the $A_{\s,1}$:
\begin{eqnarray*}
A_{\s(23),1}A_{\s,1}^{-1}& = & A_{\s,3} \\
& = & A_{\s(234),4} \\
& = & A_{\s(234)(345),5}
\\
& = & A_{\s(234)(345)(456),6} \\
&\vdots& \\
& = & A_{\s(234)(345)(456)\cdots(2n-1\,2n\,1),2} \\
& = & A_{\s(234)\cdots(2n-1\,2n\,1)(2n\,1),1}
      A_{\s(234)\cdots(2n-1\,2n\,1),1}^{-1}
\end{eqnarray*}
which may be rewritten as
\begin{equation}\label{eqA1}
A_{\s(23),1}A_{\s,1}^{-1}=
A_{\s(1\,2n-1\dots 5\,3\,2n\dots 4\,2),1}
A_{\s(2n-1\dots 3\,1)(2n\dots 4\,2),1}^{-1}
\end{equation}

Now we further reduce the set of generators by recognizing that the
$X_\s$ and $A_{\s,1}$ are redundant.

\begin{lemma}\label{ABXkllem}
For every $\s \in S_{2n}$, the generators $A_{\s,1}$, $B_\s$, and $X_\s$
depend only on $\s^{-1}(1)$ and $\s^{-1}(2)$.
\end{lemma}
\begin{proof}
Consider $\tau$ in the stabilizer of $1,2$.  Clearly $\tau$
commutes with $(1\,2)$.  In $\varphi(S_{2n})$ the stabilizer is generated
by $\set{\G_4,\G_5,\ldots\G_{2n}}$ all of which commute with
$\G_1$ and $\G_{1'}$.  So $\tau$ commutes with $\G_1$ and $\G_{1'}$ as
well. By definition $A_{\tau,1}$, $B_\tau$, and $X_\tau$ are given by
$\tau\G_1\G_{1'}\tau^{-1}$, $\tau(12)\G_{1'}\tau^{-1}$, and
$\tau(12)\G_1\tau^{-1}$ respectively.  Since $\tau$ commutes with all of
these things we have $A_{\tau,1}=A_{I,1}$, $B_\tau=B_I$, and $X_\tau=X_I$.
Suppose $\s_1^{-1}(1) = \s_2^{-1}(1)$ and $\s_1^{-1}(2) = \s_2^{-1}(2)$.
Then $\s_2=\s_1\tau$ for some $\tau$ in the stabilizer of $1,2$.
Hence $A_{\s_2,1}=A_{\s_1,1}$, $B_{\s_2}=B_{\s_1}$, and
$X_{\s_2}=X_{\s_1}$.
\end{proof}

Lemma \ref{ABXkllem} suggests a convenient trio of definitions.

\begin{defn}\label{ABXkldef}
For distinct $k,\ell\in \{1,\dots,2n\}$, $A_{k\ell}$, $B_{k\ell}$, and
$X_{k\ell}$ can be defined by
\begin{eqnarray}
A_{k\ell} & = & \s\G_1\G_{1'}\s^{-1} \label{Akl} \\
B_{k\ell} & = & \s(12)\G_{1'}\s^{-1} \label{Bkl} \\
X_{k\ell} & = & \s(12)\G_1\s^{-1} \label{Xkl}
\end{eqnarray}
where $\s \in S_{2n} = \sg{\G_2,\dots,\G_{2n}}$ is any permutation such that
$\s(k)=1$ and $\s(\ell)=2$.
\end{defn}

Now we investigate the behavior of these generators under conjugation by
elements $\s\in S_{2n}$.

\begin{proposition}
For every $\s\in S_{2n} = \sg{\G_2,\dots,\G_{2n}}$ and distinct
$k,\ell\in \{1,\dots,2n\}$, we have that
\begin{eqnarray}
\s^{-1} A_{k\ell} \s & = & A_{\s(k),\s(\ell)}\label{ASact}\\
\s^{-1} B_{k\ell} \s & = & B_{\s(k),\s(\ell)}\label{BSact}\\
\s^{-1} X_{k\ell} \s & = & X_{\s(k),\s(\ell)}\label{XSact}
\end{eqnarray}
\end{proposition}
\begin{proof}
Let $\tau\in S_{2n}$ be such that $\tau(k)=1$ and $\tau(\ell)=2$. Since
$A_{12} = \G_1\G_{1'}$ by Definition \ref{ABXkldef}, we have 
$A_{k\ell} = \tau\G_1\G_{1'}\tau^{-1}$
and $\s^{-1}A_{k \ell}\s = \s^{-1}\tau\G_1\G_{1'}\tau^{-1}\s = 
A_{\tau^{-1}\s(1),\tau^{-1}\s(2)} = A_{\s(k)\s(\ell)}$.
The same proof works for $B_{k\ell}$ and $X_{k\ell}$.
\end{proof}

From Theorem \ref{As1Xs} (with a little help from Lemma \ref{ABXkllem}),
we obtain
\begin{corollary}\label{co12}
The group $\mathcal{A}$ is generated by
$\{A_{k\ell},X_{k\ell}\}_{1 \leq k,\ell \leq 2n}$.
\end{corollary}

In terms of the $A_{ik}$ the relation \eq{eqA1} can be
written (for $\s=I$) as $A_{13}A_{12}^{-1}=A_{24}A_{34}^{-1}$.
Then conjugating by any $\s$ we get
\begin{equation}\label{Akl:rel}
A_{ij}A_{ik}^{-1}=A_{k\ell}A_{j\ell}^{-1}
\end{equation}
for any four distinct indices $i,j,k,\ell$.  With two applications of
relation \eq{Akl:rel} we see that
$A_{ij}A_{ik}^{-1}=A_{mj}A_{mk}^{-1}=A_{\ell j}A_{\ell k}^{-1}$ and so
\begin{equation}\label{Akl2:rel}
A_{ij}A_{ik}^{-1}=A_{\ell j}A_{\ell k}^{-1}
\end{equation}
for any distinct indices $i,j,k$ and $\ell\ne j,k$.  Using one more
application of \eq{Akl:rel} we can also allow $\ell=i$ in
\eq{Akl:rel} because
$A_{ij}A_{ik}^{-1}=A_{\ell j}A_{\ell k}^{-1}=A_{ki}A_{ji}^{-1}$.
In view of
\eq{Akl:rel} and \eq{Akl2:rel} and Table~\ref{Vtab-rel} we can
write a translation table for the generators $A_{\s,j}$ in terms of the
new generators $A_{k\ell}$.
\begin{tab}\label{Vtrans}
$A_{\s,j}$ in terms of $A_{k\ell}$
\begin{eqnarray}
A_{\s,1} & = & \s A_{12}\s^{-1} \\
A_{\s,3} & = & \s A_{x3}A_{x2}^{-1}\s^{-1} = \s A_{2x}A_{3x}^{-1}\s^{-1}\qquad
     where\ x\ne 2,3 \\
A_{\s,4} & = & \s A_{x4}A_{x3}^{-1}\s^{-1} = \s A_{3x}A_{4x}^{-1}\s^{-1}\qquad
     where\ x\ne 3,4 \\
\null &\vdots& \null \notag \\
A_{\s,2} & = & \s A_{x1}A_{x\,2n}^{-1}\s^{-1} = \s A_{2n\,x}A_{1x}^{-1}\s^{-1}\qquad
     where\ x\ne 1,2n
\end{eqnarray}
\end{tab}

We have reduced the generating set for $\mathcal{A}$ to
$\{X_{k\ell}, A_{k \ell} \}_{k \neq \ell}$.  Now we use the
Reidmeister--Schreier rewriting process to translate all of the
relations of $\Pitil_1$.
Using the notation of Subsection \ref{splitRMS},
$\g(\s,\G_j)=I$ and $\g(\s,\G_{j'})=A_{\s,j}^{-1}$
for $j \neq 1$.  For $j=1$, $\g(\s,\G_1)=X_\s^{-1}$, and
$\g(\s,\G_{1'})=B_\s^{-1}$. We begin by translating some of the
relations which involve $\G_1$ but not $\G_{1'}$.

$\G_1\G_1\stackrel{\tau}{\longmapsto}
\g(I,\G_1)\g(\overline{\G_1},\G_1) = X_I^{-1}X_{(12)}^{-1} =
X^{-1}_{12} X^{-1}_{21}$, so we deduce that $X_{21} = X^{-1}_{12}$
and conjugating we get
\begin{equation}\label{Xrel1}
X_{\ell k} = X^{-1}_{k \ell}.
\end{equation}

The relations $[\G_1,\G_j]$ in $\Pitil_1$ all produce the same relations
on $\sg{X_{k\ell}}$ as above.

Now we translate the triple relations (\ref{tripeq2}) for $i, j$
adjacent. We start for example with $\G_1\G_2\G_1\G_2\G_1\G_2$.
The relation $\G_1\G_2\G_1\G_2\G_1\G_2$ translates through $\tau$
to the expression
$$\g(I,\G_1) \g(\overline{\G}_1,\G_2)
\g(\overline{\G_1\G_2},\G_1) \g(\overline{\G_1\G_2\G_1},\G_2)
\g(\overline{\G_1\G_2\G_1\G_2}, \G_1)
\g(\overline{\G_1\G_2\G_1\G_2\G_1},\G_2).$$
But since $\g(\s,\G_2) = I$ we get $\g(I,\G_1)\g(\overline{\G_1\G_2},\G_1)
\g(\overline{\G_1\G_2\G_1\G_2},\G_1)$, and since
$\overline{\G_1\G_2} = (1\,2\,2n)$ we can further simplify to
$$X_I^{-1}X_{(1\,2\,2n)}^{-1}X_{(2n\,2\,1)}^{-1} =
X^{-1}_{12} X^{-1}_{2n\,1}X^{-1}_{2\,2n}.$$ Thus $X_{2\,2n}X_{2n\,1}X_{12} = 1$
and including all conjugates we have
\begin{equation}\label{Xrel2}
X_{k\ell}X_{\ell m}X_{mk} = 1.
\end{equation}

 Similarly the
relation $\G_1\G_3\G_1\G_3\G_1\G_3$ translates through $\tau$ to
the expression 
$$\g(I,\G_1) \g(\overline{\G}_1,\G_3)
\g(\overline{\G_1\G_3},\G_1) \g(\overline{\G_1\G_3\G_1},\G_3)
\g(\overline{\G_1\G_3\G_1\G_3},\G_1)
\g(\overline{\G_1\G_3\G_1\G_3\G_1},\G_3)$$ which equals
$X_I^{-1}X_{(321)}^{-1}X_{(123)}^{-1} =
X^{-1}_{12}X^{-1}_{23}X^{-1}_{31}$. Thus $X_{31}X_{23}X_{12} = 1$, and conjugating we obtain
\begin{equation}\label{Xrel3}
X_{\ell m}X_{k \ell}X_{mk} = 1.
\end{equation}

  Together the relations \eqs{Xrel1}{Xrel3} show that $\sg{X_{\ell m}}$ is
generated by the $2n{-}1$ commuting elements $X_{12},\dots,X_{1\,2n}$
so that $\sg{X_{k\ell}} \cong \Z^{2n-1}$.

Next, taking the above relations and replacing $\G_1$ with $\G_{1'}$
it is easy to see that $\G_{1'}\G_{1'}$, $[\G_{1'},\G_j]$,
$\G_{1'}\G_2\G_{1'}\G_2\G_{1'}\G_2$, and $\G_{1'}\G_3\G_{1'}\G_3\G_{1'}\G_3$
yield identical relations among the $B_{k\ell}$. 
So $\sg{B_{k\ell}} \cong \Z^{2n-1}$ as well,

\begin{eqnarray}
B_{\ell k} &=& B^{-1}_{k \ell}\label{Brel1} \\
B_{k\ell}B_{\ell m}B_{mk} &=& 1\label{Brel2} \\
B_{\ell m}B_{k \ell}B_{mk} &=& 1.\label{Brel3}
\end{eqnarray}

    We finish with the last necessary triple relations.
Note that $\overline{\G_{1'}\G_1\G_{1'}} =
\overline{\G}_1$ and $\tau(\G_{1'}\G_1\G_{1'}) =
B_I^{-1}X_{(12)}^{-1}B_I^{-1} = B_{12}^{-1}X_{21}^{-1}B_{12}^{-1}
= B_{12}^{-1}X_{12}B_{12}^{-1}$. So if we define $C_{k \ell}$ to
be $B_{k\ell}X^{-1}_{k\ell} B_{k\ell} = X_{k \ell} A^2_{k \ell}$
then the additional relations are $C_{\ell k } = C^{-1}_{k \ell}$,
$C_{k\ell}C_{\ell m}C_{mk} = 1$, and $C_{\ell m}C_{k \ell}C_{mk} = 1$.
By the arguments above, the $\{C_{k \ell}\}$ generate another
copy $\Z^{2n-1} \subset \mathcal{A}$.  In fact for each exponent $n$
the elements $X_{k \ell}A^n_{k \ell}$ generate a
subgroup  isomorphic to $\Z^{2n-1}$.

     The relations computed thus far turn out to be all of the relations
in $\pi_1(X_{\Gal}^{\Aff})$.  Computations showing that the remaining
relations translated from $\Pitil_1$ are consequences of the relations
above are identical to computations in \cite{AGTV} and are omitted here.
We have therefore proven the following theorem.

\begin{theorem}\label{th100}
The fundamental group $\pi_1(X_{\Gal}^{\Aff})$ is generated by
elements $\{X_{ij}, A_{ij}\}$ with the relations
\begin{eqnarray}
X_{ji}A_{ji}^n & =& (X_{ij}A_{ij}^n)^{-1},\\
(X_{ij}A_{ij}^n)(X_{jk}A_{jk}^n)(X_{ki}A_{ki}^n)& =& 1,\\
(X_{jk}A_{jk}^n)(X_{ij}A_{ij}^n)(X_{ki}A_{ki}^n)& =& 1,\\
A_{ij}A_{ik}^{-1} & = & A_{k\ell}A_{j\ell}^{-1}
\end{eqnarray}
for every $n \in \Z$ and distinct $i, j, k, \ell$.
\end{theorem}

Before adding the projective relation to compute $\pi_1(X_{\Gal})$
we prove a useful lemma, showing that some of the $A_{k\ell}$
commute in $\pi_1(X_{\Gal}^{\Aff})$.  We shall frequently use the fact that
$X_{k\ell}X_{\ell m}=X_{\ell m}X_{k\ell}=X_{km}$ which is a consequence
of \eqs{Xrel1}{Xrel3}.  $B_{k\ell}$ and $C_{k\ell}$ satisfy this as
well.

\begin{lemma}\label{Acomlem}
In $\pi_1(X_{\Gal}^{\Aff})$ we have $[A_{ij},A_{ik}]=1$ and
$[A_{ji},A_{ki}]=1$ for distinct $i,j,k$.
\end{lemma}

\begin{proof}
Starting with $1=C_{ki}C_{jk}C_{ij}$ and use the definition of $C_{ij}$
to rewrite it as
\begin{align*}  
1&=(B_{ki}X_{ik}B_{ki})(B_{jk}X_{kj}B_{jk})(B_{ij}X_{ji}B_{ij})
=B_{ki}X_{ik}B_{ji}X_{kj}B_{ik}X_{ji}B_{ij}\\
&=B_{ki}X_{ik})(B_{ji}X_{ij})(X_{ki}B_{ik})(X_{ji}B_{ij})=
A_{ik}^{-1}A_{ij}^{-1}A_{ik}A_{ij}.  
\end{align*}
Thus the commutator
$[A_{ij},A_{ik}]=1$.  The relation (\ref{Akl:rel}) can be used to show
that the commutator $[A_{ji},A_{ki}]=1$ as well.
\end{proof}

\subsection{The projective relation}\label{sec:9}

     To complete the computation of $\pi_1(X_{\Gal})$ we need only to add
the projective relation
$$\G_1\G_{1'}\G_2\G_{2'}\cdots\G_{2n}\G_{2n'} = 1.$$
This relation translates in $\mathcal{A}$ as the product
$P = A_{I,1}A_{I,2}\cdots A_{I,2n}$.
We must translate the $A_{I,j}$ in terms of the $A_{k
\ell}$, using Table \ref{Vtrans}.  $P$ translates to
$$A_{12}(A_{31}A_{3\,2n}^{-1}) (A_{21}A_{31}^{-1})(A_{31}A_{41}^{-1})\cdots
(A_{2n-2\,1}A_{2n-1\,1}^{-1})(A_{2n-1\,1}A_{2n\,1}^{-1}).$$
All but five terms cancel in the expression above, leaving
$A_{12}A_{31}A_{3\,2n}^{-1}A_{21}A_{2n\,1}^{-1}$.  Using Equation
\eq{Akl:rel}, we get $A_{12}A_{31}A_{3\,2n}^{-1}A_{3\,2n}A_{32}^{-1}$.
Thus the projective relation may be written as
$A_{12}A_{31}A_{32}^{-1}=1$ or equivalently $A_{32}=A_{12}A_{31}$.
Conjugating, this becomes
\begin{equation}\label{Akl:rel2}
A_{ij} = A_{kj}A_{ik}.
\end{equation}
Substituting back into \eq{Akl:rel}, writing $A_{ij} =
A_{kj}A_{ik}$ and $A_{k\ell} = A_{j\ell}A_{kj}$, we obtain
\begin{equation}\label{Akl:rel3}
A_{kj} A_{j \ell} = A_{j\ell}A_{kj}.
\end{equation}

\begin{lemma}
The subgroup $\sg{A_{k\ell}}$ of $\pi_1(X_{\Gal})$ is
commutative of rank of at most $2n-1$.
\end{lemma}
\begin{proof}
We will compute the centralizer of $A_{ij}$ for fixed $i,j$. Let
$i,j,k,\ell$ be four distinct indices.  We already know from Lemma
\ref{Acomlem} that $A_{ij}$ commutes with $A_{ik}$ and $A_{\ell j}$.
By equation (\ref{Akl:rel3}) it also commutes with $A_{ki}$ and
$A_{j\ell}$.  Now equation (\ref{Akl:rel2}) allows us to write
$A_{k\ell}=A_{i\ell}A_{\ell j}$, both of which commute with $A_{ij}$,
so $\sg{A_{k\ell}}$ is commutative.

Now, since $A_{j k} A_{ij} = A_{ik}$ and $A_{ik}A_{ji} = A_{jk}$, we have
$A_{jk} A_{ij} A_{ji} = A_{ik} A_{ji} = A_{jk}$, so that $A_{ji} =
A_{ij}^{-1}$, the group is generated by the $A_{1k}$ ($k = 2,\dots,2n$), and
the rank is at most $2n-1$.
\end{proof}

We see that $\pi_1(X_{\Gal}) = \sg{A_{ij},X_{ij}}$ with the two
  subgroups $\sg{A_{ij}},\sg{X_{ij}}$ each isomorphic to  $\Z^{2n-1}$.
The only question left is how these two subgroups interact.

\begin{lemma}\label{lm101}
In $\pi_1(X_{\Gal})$ the $A_{ij}$ and $X_{k\ell}$ commute.
\end{lemma}

\begin{proof}
We need only consider the commutators of $A_{12}$ and $X_{ij}$
since all others are merely conjugates of these. First consider
the commutator $[X_{12},A_{12}]$. Since
$X_{12}=(12)\G_1$ and $A_{12}=\G_1\G_{1'}$ the commutator
$X_{12}A_{12}X_{12}^{-1}A_{12}^{-1}$ becomes
$(12)\G_1(\G_1\G_{1'})\G_1(12)A_{12}^{-1}=(12)\G_{1'}\G_1(12)A_{12}^{-1}=
A_{21}^{-1}A_{12}^{-1}=A_{12}A_{12}^{-1}=1$.  So $X_{12}$ and
$A_{12}$ commute.

Next consider the commutator $X_{13}A_{12}X_{13}^{-1}A_{12}^{-1}$.
By definition we have that 
$X_{13}=(23)X_{12}(23)=(23)(12)\G_1(23)=(123)\G_1\G_3$.
We are using the fact that $\G_3=(23)$ in the image $\varphi(S_{2n})$.
Thus the commutator can be written as
$(123)\G_1\G_3(\G_1\G_{1'})\G_3\G_1(321)A_{12}^{-1}$. We use the
triple relations $$\Trip{\Gamma_{(1)}}{\Gamma_{(3)}}$$ to rewrite it as
\begin{align*}  
(123)\G_3\G_1\G_3\G_{1'}\G_3\G_1(321)A_{12}^{-1}&=
(123)\G_3\G_1\G_{1'}\G_3\G_{1'}\G_1(321)A_{12}^{-1}\\
&=(123)(23)\G_1\G_{1'}(23)(321)(123)\G_{1'}\G_1(321)A_{12}^{-1}
\end{align*}
which is equal to $A_{(123)(23),1}A_{(123),1}^{-1}A_{12}^{-1}=
A_{(13),1}A_{(123),1}^{-1}A_{12}^{-1}=A_{32}A_{31}^{-1}A_{12}^{-1}=
A_{32}A_{32}^{-1}=1$, 
proving that $X_{13}$ and $A_{12}$ commute.
Conjugating by $(3j)$ we see that $X_{1j}$ commutes with $A_{12}$
and since $X_{ij}=X_{1i}^{-1}X_{1j}$ we see that every $X_{ij}$
commutes with $A_{12}$.
\end{proof}

\begin{theorem}\label{main}
The fundamental group $\pi_1(X_{\Gal})\cong\Z^{4n-2}$.
\end{theorem}
\begin{proof}
$\pi_1(X_{\Gal})$ is generated by $A_{1j}$ and $X_{1j}$ which all
commute.  Hence the group they generate is $\Z^{4n-2}$.
\end{proof}

\subsection*{Acknowledgements}
\addcontentsline{toc}{section}{Acknowledgements}

The work is partially supported by the Golda Meir Fellowship,
Mathematics Institute, Hebrew university, Jerusalem. The first author
wishes to thank the Mathematics Institute, Hebrew university,
Jerusalem and her present host Hershel Farkas.

\Addresses


\begin{thebibliography}

\bibitem{AGTV} Amram, M., Goldberg, D., Teicher, M., Vishne, U., {\em
The fundamental group of a Galois cover of $\C\P^1\times T$},
\href{http://www.maths.warwick.ac.uk/agt/AGTVol2/agt-2-20.abs.html}{Algebr.
Geom. Topol. {\bf 2}, (2002) 403-432.} \MR{1917060}

\bibitem{Mo} Moishezon, B., {\em Algebraic surfaces and the arithmetic
of braids, II}, Combinatorial methods in topology and algebraic
geometry (Rochester, N.Y., 1982), Contemp. Math. {\bf 44}, (1985)
311-344. \MR{0813122}

\bibitem{MoTe2}
 Moishezon, B., Teicher, M., {\it Braid group
   technique in complex geometry I, Line arrangements in $\C\P^2$},
   Braids (Santa Cruz, CA, 1986),
   Contemp. Math. {\bf 78}, (1988)  425-555. \MR{0975093}

\bibitem{MoTe4} Moishezon, B., Teicher, M., {\it Braid group technique
 in complex geometry IV: Braid monodromy of the branch curve $S_3$ of
 $V_3 \rightarrow \C\P^2$ and application to $\pi_1(\C\P^2 - S_3,
 \ast)$,} Classification of algebraic varieties (L'Aquila, 1992),
 Contemp. Math. {\bf 162}, (1993) 332-358. \MR{1272707}

\bibitem{vK}
van Kampen, E.R., {\em On the fundamental group of an algebraic curve},
Amer. J. Math. {\bf 55}, (1933) 255-260.

\end{thebibliography}
\end{document}